\documentclass[11pt]{amsart}
\usepackage{amsmath,amssymb,amsthm,amscd}
\usepackage{palatino, mathpazo}
\usepackage{amsfonts}
\usepackage[mathscr]{eucal}
\usepackage[all]{xy}

\newtheorem{theorem}{Theorem}[section]
\newtheorem{lemma}[theorem]{Lemma}
\newtheorem{proposition}[theorem]{Proposition}
\newtheorem{corollary}[theorem]{Corollary}

\theoremstyle{remark}
\newtheorem{remark}[theorem]{Remark}
\newtheorem{example}[theorem]{Example}
\numberwithin{equation}{section}

\newcommand{\T}{\mathscr{F}}

\def\dim{\mathop{\rm dim}\nolimits}

\def\cit{{\mathbb C}}

\def\zit{{\mathbb Z}}
\def\pit{{\mathbb P}}
\def\0{{\mathcal O}}

\def\X{{\mathfrak X}}

\numberwithin{equation}{section}

\renewcommand{\O}{{\mathscr O}}

\newcommand{\F}{{\mathfrak F}}

\newcommand{\E}{{\mathscr E}}

\newcommand{\A}{{\mathbb A}}

\newcommand{\Hom}{\operatorname{Hom}}
\newcommand{\codim}{\operatorname{codim}}

\newcommand{\Id}{\operatorname{Id}}
\newcommand{\End}{\operatorname{End}}

\newcommand{\im}{\operatorname{Im}}
\def\S{Section~}

\title[stratified Mukai flops]
{Motivic and quantum invariance\\ under\\ stratified Mukai flops}

\author{Baohua~Fu}
\address{B.\ Fu; C.N.R.S., Labo.\ J.\ Leray, Facult¡äe des sciences,
Universit\'e de NANTES, 2, Rue de la Houssini\`ere, BP 92208,
F-44322 Nantes Cedex 03 France.} \email{fu@math.univ-nantes.fr}

\author{Chin-Lung~Wang}
\address{C.-L.\ Wang; Dept.\ of Math., National Central
University, Jhongli, Taiwan. National Center for Theoretic
Sciences, Hsinchu, Taiwan.} \email{dragon@math.ncu.edu.tw}

\begin{document}
\maketitle

\begin{abstract}
For stratified Mukai flops of type $A_{n,k}, D_{2k+1}$ and
$E_{6,I}$, it is shown the fiber product induces isomorphisms on
Chow motives.

In contrast to (standard) Mukai flops, the cup product is
generally not preserved. For $A_{n, 2}$, $D_5$ and $E_{6, I}$
flops, quantum corrections are found through
degeneration/deformation to ordinary flops.
\end{abstract}
\section{Introduction}

\subsection{Backgroumd}

Two smooth projective varieties over $\mathbb{C}$ are \emph{$K$
equivalent} if there are birational morphisms $\phi: Y \to X$ and
$\phi': Y \to X'$ such that $\phi^*K_X = \phi'^* K_{X'}$. This
basic equivalence relation had caught considerable attention in
recent years through its appearance in minimal model theory,
crepant resolutions, as well as other related fields.

The conjectural behavior of $K$ equivalence has been formulated in
\cite{W}. A \emph{canonical correspondence} $\T \in A^*(X \times
X')$ should exist and gives an isomorphism of Chow motives $[X]
\cong [X']$. Under $\T$, $X$ and $X'$ should have isomorphic
\emph{B-models} (complex moduli with Hodge theory on it) as well
as \emph{A-models} (quantum cohomology ring up to analytic
continuations over the extended K\"ahler moduli space).

Basic examples of $K$ equivalence are \emph{flops} (with
exceptional loci $Z, Z', S$)
$$
\xymatrix{(X, Z) \ar[rd]^{\psi} \ar@{.>}[rr]^f & & (X', Z')
\ar[ld]_{\psi'}\\ & (\bar X, S) &}.
$$
Among them the \emph{ordinary flops} had been studied in
\cite{LLW} where the equivalence of motives and A-models was
proved. In that case $\T$ is the graph closure $\T_0 :=
\bar\Gamma_f$. In general, $\T$ must contain \emph{degenerate
correspondences}.

The typical examples are \emph{Mukai flops}. They had been
extensively studied in the literature in hyper-K\"ahler geometry.
Over a general base $S$, they had also been studied in \cite{LLW},
where the invariance of Gromov-Witten theory was proved. In that
case $\T = X \times_{\bar X} X' = \T_0 + \T_1$ with $\T_1 = Z
\times_S Z'$. We expect that for flops $\T$ should be basically $X
\times_{\bar X} X'$.

To understand the general picture we are led to study flops with
$\T$ consisting of many components. The \emph{stratified Mukai
flops} provide such examples. They appear naturally in the study
of symplectic resolutions \cite{Fu} \cite{Na} and they should play
important roles in higher dimensional birational geometry. For
hyper-K\"ahler manifolds, see for example \cite{Mar}.

In this paper, we study general stratified Mukai flops without any
assumptions on the global structure of $X$ and $X'$. By way of it,
we hope to develop tools with perspective on future studies.

\subsection{Stratified Mukai flops}

Fix two natural numbers $n, k$ such that $2k < n+1$. Consider two
smooth projective varieties $X$ and $X'$. Let $F_k \subset F_{k-1}
\subset \cdots F_1 \subset X$ and $F'_k \subset F'_{k-1} \subset
\cdots F'_1 \subset X'$ be two collections of closed subvarieties.
Assume that there exist two birational morphisms $X
\xrightarrow{\psi} \bar X \xleftarrow{\psi'} X'$. The induced
birational map $f: X \dasharrow X'$ is called {\em a (stratified)
Mukai flop of type $A_{n,k}$} over $\bar{X}$ if the following
conditions are satisfied:
\begin{itemize}
\item[(i)] The map $f$ induces an isomorphism $X \setminus F_1
\mathop{\longrightarrow}\limits^\sim X' \setminus F'_1$;

\item[(ii)] $\psi(F_j) = \psi'(F'_j) =:S_j$ for $1 \leq j \leq k$;

\item[(iii)] $S_k$  is smooth and there exists a vector bundle $V$
of rank $n + 1$ over it such that $F_k$ is isomorphic to the
relative Grassmanian $ G_{S_k}(k, V)$ of $k$-planes over $S_k$ and
the restriction $ \psi|_{F_k}: F_k \to S_k$ is the natural
projection. Furthermore, the normal bundle $N_{F_k/X}$ is
isomorphic to the relative cotangent bundle $T^*_{F_k/S_k}$. The
analogue property holds for $F'_k$ and $\psi'$ with $V$ replaced
by its dual $V^*$;

\item[(iv)] If $k=1$, we require that $f$ is a usual Mukai flop along
$F_k$. When $k \geq 2$, let $Y$ (resp. $Y'$, $\bar{Y}$) be the
blow-up of $X$ (resp. $X'$, $\bar{X}$) along $F_k$ (resp. $F'_k$,
$S_k$). By the universal property of the blow-ups, we obtain
morphisms $Y \to \bar{Y} \leftarrow Y'$. The proper transforms of
$F_j$, $F'_j$ give collections of subvarieties on $Y$, $Y'$.  We
require that the birational map $Y \dasharrow Y'$ is a Mukai flop
of type $A_{n-2, k-1}$.
\end{itemize}

We define a {\em Mukai flop of type $D_{2k+1}$} in a similar way
with the following changes: (1) one requires that $S_k$ is simply
connected; (2) the vector bundle $V$ is of rank $4k+2$ with a
fiber-wise non-degenerate symmetric 2-form. Then the relative
Grassmanians of $k$-dimensional isotropic subspaces of $V$ over
$S_k$ has two components $G_{iso}^+$ and $G_{iso}^-$. We require
that $F_k$ (resp. $F'_k$) is isomorphic to $G_{iso}^+$ (resp.
$G_{iso}^-$); (3) when $k=1$, $f$ is a usual Mukai flop.

Similarly one can define a {\em Mukai flop of type $E_{6, I}$} by
taking $k=2$ with $V$ being an $E_6$-vector bundle of rank 27 over
$S_2$ and $F_2$ is the relative $E_6/P_1$-bundle over $S_2$ in
$\pit(V)$. The dual variety $F'_2$ is given by the relative
$E_6/P_6$-bundle in $\pit(V^*)$, where $P_1, P_6$ are maximal
standard parabolic subgroup in $E_6$ corresponding to the simple
roots $\alpha_1, \alpha_6$ respectively. By \cite{CF}, when we
blow up the smallest strata of the flop, we obtain a usual Mukai
flop.

\subsection{Main results}

Our main objective of this work is to prove the following
theorems.

\begin{theorem}\label{main}
Let $f: X \dasharrow X'$ be a Mukai flop of type $A_{n, k}$,
$D_{2k+1}$ or $E_{6, I}$ over $\bar{X}$. Let $\T$ be the
correspondence $X \times_{\bar{X}} X'$. Then $X$ and $X'$ have
isomorphic Chow motives under $\T$. Moreover $\T$ preserves the
Poincar\'e pairing of cohomology.
\end{theorem}

Note that the flops of type $A_{n, 1}$ and $D_3$, i.e.~$k = 1$,
are the usual Mukai flops, and in these cases the theorem has been
proven in \cite{LLW}. Our proof uses an induction on $k$ (for all
$n$) via (iv). We shall give details of the proof for $A_{n,k}$
flops, while omitting the proof of the other two types, since the
argument is essentially the same.

For $k = 1$ (i.e.~the usual Mukai flops), the cohomology ring as
well as the Gromov-Witten theory are also invariant under $\T$
\cite{LLW}. However, the general situation is more subtle:

\begin{theorem} \label{main-2}
When $k \ge 2$, the cup product is generally not preserved under
$\T$. For $A_{n, 2}$, $D_5$ and $E_{6, I}$ flops the defect is
corrected by the genus zero Gromov-Witten invariants attached to
the extremal ray, up to analytic continuations.
\end{theorem}

While Theorem \ref{main} is as expected, Theorem \ref{main-2} is
somehow surprising, since stratified Mukai flops are in some sense
locally (holomorphically) symplectic and it is somehow expected
that there is no quantum corrections for flops of these types.
Indeed stratified Mukai flops among hyper-K\"ahler manifolds can
always be deformed into isomorphisms \cite{Huy} hence there is no
quantum correction. As it turns out, the key point is that for the
projective local models of general stratified Mukai flops, in
contrast to the case $k = 1$, we cannot deform them into
isomorphisms!

\subsection{Outline of the contents}

In \S2, the existence of $A_{n, k}$ flops in the projective
category is proved via the cone theorem. In \S3, a general
criterion on equivalence of Chow motives \emph{via graph closure}
is established for \emph{strictly semi-small} flops. While a given
flop may not be so, generic deformations of it may sometimes do.
When this works, we then restrict the graph closure of the one
parameter deformation back to the central fiber to get the
correspondence, which is necessarily the fiber product.

It is thus crucial to study deformations of flops. Global
deformations are usually obstructed, so instead we study in \S4
the \emph{deformations of projective local models} of $A_{n, k}$
flops. While open local models can be deformed into isomorphisms,
the projective local models cannot be deformed into isomorphisms
in general but only be deformed into certain $A^*_{n - 2, k - 1}$
flops.  These flops, which we called \emph{stratified ordinary
flops}, do not seem to be studied before in the literature.
Nevertheless this deformation is good enough for applying the
equivalence criterion of motives.

To handle global situations, we consider \emph{degenerations to
the normal cone} to reduce problems on $A_{n, k}$ flops to
problems on $A_{n - 2, k - 1}$ flops and on local models of $A_{n,
k}$ flops. This is carried out for correspondences in \S5. This
makes inductive argument work since local models are already well
handled. We also carry out this for \emph{cup product} by proving
an \emph{orthogonal decomposition} under degenerations to the
normal cone. This in particular applies to the Poincar\'e pairing
and completes the proof of Theorem \ref{main}.

In section 6 we prove Theorem \ref{main-2} for $A_{n, 2}$ flops.
We apply the \emph{degeneration formula} for GW invariants
\cite{LR} \cite{Li} to splits the \emph{absolute} GW invariants
into the \emph{relative} ones. After degenerations, the flop is
split into \emph{two simpler flops}, one is a Mukai flop and
another one can be deformed into ordinary $\mathbb{P}^{n - 2}$
flop. It turns out that each \emph{GW invariant attached to the
extremal ray} must go to one of these two factors completely. For
the former the extremal invariants indeed vanish. For the latter
we use a recent result on ordinary flop with general base
\cite{LLW2} to achieve the quantum corrections up to analytic
continuations. This then completes the proof.

At the end we compare Theorem \ref{main-2} with the hyper-K\"ahler
case, where the ring structure is preserved and there are no
non-trivial Gromov-Witten invariants. When $f$ is not standard
Mukai, all these may fail without the global hyper-K\"ahler
condition. A careful comparison of the degeneration analysis in
this case with the local model case leads to some new topological
constraint on hyper-K\"ahler manifolds (c.f.~Proposition
\ref{main-3}).

\subsection{Acknowledgements}

B. Fu is grateful to the Department of Mathematics, National
Central University (Jhongli, Taiwan) for providing excellent
environment which makes the collaboration possible.

C. L. Wang would like to thank the MATHPYL program of the
F\'ed\'eration de Math\'ematiques des Pays de Loire for the
invitation to Nantes.

\section{Existence of (twisted) $A_{n, k}$ flops}

Given $k, n \in \mathbb{N}$, $2k < n + 1$, a flopping contraction
$\psi: (X, F) \to (\bar X, S)$ is of type $A_{n, k}$ if it admits
the following inductive structure: There is a filtration $F = F_1
\supset \cdots \supset F_k$ with induced filtration $S = S_1
\supset \cdots \supset S_k$, $S_j := \psi(F_j)$ such that
$\psi_{S_k}: F_k \cong G_{S_k}(k, V) \to S_k$ is a $G(k, n + 1)$
bundle for some vector bundle $V \to S_k$ of rank $n + 1$ with
$$
N_{F_k/X} \cong T^*_{F_k/S_k} \otimes \psi_{S_k}^* L_k
\quad\mbox{for some}\quad L_k \in {\rm Pic}\,S_k.
$$
Moreover, the blow-up maps $\phi$, $\bar\phi$ fit into a cartesian
diagram
$$
\xymatrix{Y = {\rm Bl}_{F_k} X \supset E \ar[d]_\phi
\ar[rd]^{\bar\psi} & \\
X \ar[rd]_\psi & \bar Y = {\rm Bl}_{S_k} \bar X \supset \bar E
\ar[d]^{\bar\phi} \\
& \bar X}
$$
such that the induced contractions $\bar \psi: (Y, \tilde F) \to
(\bar Y, \tilde S)$ with filtrations $\tilde F = \tilde F_1
\supset \cdots \supset \tilde F_{k - 1}$, $\tilde F_j :=
\phi_*^{-1}(F_j)$, $\tilde S = \tilde S_1 \supset \cdots \supset
\tilde S_{k - 1}$, $\tilde S_j = \bar\phi_*^{-1}(S_j)$, $1 \le j
\le k - 1$ is of type $A_{n - 2, k - 1}$. Here we use the
convention that an $A_{n, 0}$ contraction is an isomorphism. By
definition, $A_{n, 1}$ contractions are twisted Mukai
contractions.

The main results of this paper are all concerned with the
(untwisted) stratified Mukai flops, namely $L_k \cong \O_{S_k}$.
The starting basic existence theorem of flops does however hold
for the twisted case too.

\begin{proposition}
Given any $A_{n, k}$ contraction $\psi$, the corresponding $A_{n,
k}$ flop
$$
\xymatrix{(X, F) \ar[rd]_\psi \ar@{.>}[rr]^f &
& (X', F') \ar[ld]^{\psi'} \\
& (\bar X, S) & }
$$
exists with $\psi'$ being an $A_{n, k}$ contraction.
\end{proposition}

\begin{proof}
We construct the flop by induction on $k$. The case $k = 1$ has
been done in \cite{LLW}, section 6, so we let $k \ge 2$ and $n + 1
> 2k$. By induction we have a diagram
$$
\xymatrix{Y \ar[d]_\phi \ar[rd]^{\bar\psi} \ar@{.>}[rr]^g & & Y'
\ar[ld]_{\bar \psi'}
\ar@{.>}[d]_?^{\phi'}\\
X \ar[rd]_\psi & \bar Y \ar[d]^{\bar\phi} & X'\\
& \bar X &}
$$
where $g: (Y, \tilde F) \dasharrow (Y', \tilde F')$ is an $A_{n -
2, k - 1}$ flop and $\bar\psi': (Y', \tilde F') \to (\bar Y,
\tilde S)$ is an $A_{n - 2, k - 1}$ contraction.

Let $C \subset E$ be a $\phi$-exceptional curve and $C' = g_* C$
be its proper transform in $E' = g_* E$. We shall construct a
blow-down map $\phi': Y' \to X'$ for $C'$. Let $\gamma$ (resp.\
$\gamma'$) be the flopping curve for $\bar\psi$ (resp.\
$\bar\psi'$).

Since the Poincar\'e pairing is trivially preserved by the graph
correspondence $\T_0$ of $g$ in the divisor/curve level, and $\T_0
C = C' + a\gamma'$ for some $a \in \mathbb{N}$ (in fact $a = 1$),
we compute
$$
(K_{Y'}.C') = (K_{Y'}.\T_0 C) = (K_Y.C) < 0.
$$
To show that $C'$ is an Mori (negative) extremal curve, it is thus
sufficient to find a supporting divisor for it.

Let $\bar L$ be a supporting divisor for $\bar C = \bar\psi(C)$ in
$\bar Y$. Then $\bar\psi'^* \bar L$ is a supporting divisor for
the extremal face spanned by $C'$ and $\gamma'$. The idea is to
perturb it to make it positive along $\gamma'$ while keeping it
vanishing along $C'$.

Let $D$ be a supporting divisor for $C$ in $Y$ with $\lambda :=
(D.\gamma) > 0$. Let $D' = g_* D = \T_0 D$. Since $\T_0 \gamma = -
\gamma'$, we compute
$$
(D'.\gamma') = -(D, \gamma) = -\lambda < 0,
$$
$$
(D'. C') = (D'.\T_0 C) - a(D', \gamma') = (D.C) + a\lambda =
a\lambda.
$$

Let $H'$ be a supporting divisor for $\gamma'$ in $Y'$ with $c' :=
(H'.C') > 0$. Then
$$
W := a\lambda H' - c'D'
$$
has the property that $(W.\gamma') > 0$ and $(W.C') = 0$.

Now for $m$ large enough, the perturbation
$$
L' := m\bar\psi'^* \bar L + W
$$
is a supporting divisor for $C'$. Indeed, $L'$ takes the same
values as $W$ on $\gamma'$ and $C'$, while $(L'.\beta') > 0$ for
other curve classes $\beta'$ in $Y'$. That is, $L'$ is big and nef
which vanishes precisely on the ray $\mathbb{Z}^+[C']$.

By the (relative) cone theorem applying to
$\bar\phi\circ\bar\psi': Y' \to \bar X$, we complete the diagram
and achieve the flop $f: X \dasharrow X'$:
$$
\xymatrix{Y \ar[d]_\phi \ar[rd]^{\bar\psi} \ar@{.>}[rr]^g & & Y'
\ar[ld]_{\bar \psi'}
\ar[d]^{\phi'}\\
X \ar[rd]_\psi \ar@{.>}[rr]|{\qquad}
& \bar Y \ar[d]^{\bar\phi} & X' \ar[ld]^{\psi'}\\
& \bar X &}
$$

It remains to show that the contraction $\psi': X' \to \bar X$ is
of type $A_{n, k}$. By construction, it amounts to analyze the
local structure of $F_k' := \phi'(E')$. Since the flop $f$ is
unique and local with respect to $\bar X$, it is enough to
determine its structure in a neighborhood of $S_k$. This can be
achieved by explicit constructions.

Suppose that $F_k = G_{S_k}(k, V)$. We consider the pair of spaces
$(\tilde X', \tilde F_k')$ defined by duality. Namely $\tilde F_k'
:= G_{S_k}(k, V^*)$ and
$$
\mbox{$\tilde X'$ is the total space of} \quad T^*_{\tilde
F_k/S_k} \otimes \psi_{S_k}^* L_k.
$$
It is well-known that, in a neighborhood of $S_k$, $X \dasharrow
X'_0$ is an $A_{n, k}$ flop. Thus the local structure of $(X',
F_k')$ must agree with $(\tilde X', \tilde F_k')$. The proof is
complete.
\end{proof}

\begin{remark}
In the definition of $A_{n, k}$ contractions, the restriction to
exceptional divisors $\bar\psi|_E: (E, \tilde F|_E) \to (\bar E,
\tilde S|_{\bar E})$ is also an $A_{n - 2, k - 1}$ contraction.
Moreover, in the proposition the restriction
$$
\xymatrix{(E, \tilde F|_E) \ar@{.>}[rr]^{g|_E}
\ar[rd]^{\bar\psi|_E} & & (E', \tilde F'|_{E'}
\ar[ld]_{\bar\psi'|_{E'}}) \\
& (\bar E, \tilde S|_{\bar E})}
$$
is also an $A_{n - 2, k - 1}$ flop.
\end{remark}

\section{Equivalence criteria of motives}
Let $X \xrightarrow{\psi} \bar X \xleftarrow{\psi'} X'$ be two
projective resolutions of a quasi-projective normal variety $\bar
X$, and $f: X \dasharrow X'$ the induced birational map. Consider
the graph closure $\Gamma$ of $f$ and $X \xleftarrow{\phi} \Gamma
\xrightarrow{\phi'} X'$ the two graph projections. Then we obtain
a morphism between Chow groups:
$$
\T := \phi'_* \phi^*: A^*(X) \to A^*(X').
$$
For any $i$, we will consider the closed subvariety
$$
E_i= \{ x \in X \mid \dim_x \psi^{-1}(\psi(x)) \geq i \}.
$$
In a similar way we define the subvariety  $E'_i$ on $X'$. By
Zariski's main theorem, $\psi$ is an isomorphism over $X \setminus
E_1$, thus $\psi(E_1) = \psi'(E'_1) = \bar X_{sing}$.

The following criterion generalizes the one for ordinary flops in
\cite{LLW}:

\begin{proposition} \label{equi}
If for any irreducible component $D$, $D'$ of $E_i$ and $E'_i$
respectively, we have
$$
2 \codim D > \codim \psi(D), \text{and}\ 2 \codim D' > \codim
\psi'(D'),
$$
then $\T$ is an isomorphism on Chow groups which preserves the
Poincar\'e pairing on cohomology groups.

Moreover, the correspondence $[\Gamma]$ induces an isomorphism
between Chow motives: $[X] \simeq [X']$.
\end{proposition}

\begin{proof}
For any smooth $T$, $f\times {\rm id}_T: X \times T
\dashrightarrow X' \times T$ is also a birational map with the
same condition. Thus by the identity principle we only need to
prove the equivalence of Chow groups under $\T$.

For any $\alpha \in A_k(X)$, up to replacing by an equivalent
cycle, we may assume that $\alpha$ intersects $E:=\sum_{i \geq 1}
E_i$ properly. Then we have $\T \alpha = \alpha'$, where $\alpha'$
is the proper transform of $\alpha$ under $f$. If we denote by
$\tilde{\alpha}$ the proper transform of $\alpha$ in
$A^*(\Gamma)$, then we have
$$
\phi'^* \alpha' = \tilde{\alpha} + \sum_C a_C F_C,
$$
where $F_C$ are some irreducible $k$-dimensional subvariety in
$\Gamma$ and $a_C \in \zit$.

For any $C$, note that $\phi'(F_C)$ is contained in the support of
$\alpha' \cap E'_1$. As $\psi'(\alpha' \cap E'_1) = \psi(\alpha)
\cap \bar X_{sing}  = \psi(\alpha \cap E_1)$, $F_C$ is contained
in $\phi^{-1} \psi^{-1}(B_C)$, where $B_C := \psi \phi(F_C)
\subset \psi(\alpha \cap E_1)$. Take the largest $i$  such that
there exists an irreducible component $D$ of $E_i$ with  $B_C
\subset \psi(\alpha \cap D)$. For a general point $s \in B_C$, we
denote by $F_{C,s}$ its fiber by the map $\psi \circ \phi$. Then
we have
$$
\dim F_{C,s} \geq \dim F_C - \dim B_C \geq \dim F_C -
\dim (\alpha \cap D) = \codim D.
$$

By our assumption, we have $\codim D > \dim D - \dim \psi(D)$, the
latter being the dimension of a general fiber of $\psi^{-1}(B_C)
\to B_C$. Thus the general fiber of the map $\phi|_{F_C}$ has
positive dimension, which gives that $\phi_* (F_C) = 0$. This
gives that $\T' \circ \T = {\rm Id}$, where $\T' = \phi_*
\phi'^*$. A similar argument then shows that $\T \circ \T' = {\rm
Id}$, thus $\T$ and $\T'$ are isomorphisms.

Since $F_C$ has positive fiber dimension in both $\phi$ and
$\phi'$ directions, the statement on Poincar\'e pairing follows
easily as in \cite{LLW}, Corollary 2.3.
\end{proof}

Now consider two (holomorphic) symplectic resolutions: $X
\xrightarrow{\psi} \bar X \xleftarrow{\psi'} X'$. A conjecture in
\cite{FN} asserts that $\psi$ and $\psi'$ are deformation
equivalent, i.e.\ there exist deformations of $\psi$ and $\psi'$
over $\cit$: $\X \xrightarrow{\Psi} \bar \X \xleftarrow{\Psi'}
\X'$, such that for any $t \neq 0$, the morphisms $\Psi_t,
\Psi'_t$ are isomorphisms. This conjecture has been proved in
various situations, such as nilpotent orbit closures of classical
type \cite{FN} \cite{Na}, or when $W$ is projective \cite{Na1}.

Assume this conjecture and consider the birational map $\F: \X
\dasharrow \X'$. Recall that every symplectic resolution is
automatically semi-small by the work of Kaledin \cite{Ka} and
Namikawa \cite{Na1}. We obtain that the deformed resolutions
$\Psi$ and $\Psi'$ satisfy the condition of the precedent
proposition. As a consequence, we obtain:

\begin{theorem}
Consider two symplectic resolutions $X \xrightarrow{\psi} \bar X
\xleftarrow{\psi'} X'$. Suppose that they are deformation
equivalent (say given by $\F: \X \dasharrow \X'$). If we denote by
$\Gamma$ the graph of $\F$ and $\Gamma_0$ its central fiber. Then
the correspondence $[\Gamma_0]$ induces an isomorphism of motives
$[X] \simeq [X']$ which preserves also the Poincar\'e pairing.
\end{theorem}

\section{Deformations of local models}

{}From now on all the stratified Mukai flops are untwisted.

\subsection{Deformations of open local models}
Let $S$ be a smooth variety and $V \to S$ a vector bundle of rank
$n+1$. The relative Grassmanian bundle of $k$-planes in $V$ is
denoted by $\psi: F:= G_S(k,V) \to S$. Let $T$ be the universal
sub-bundle of rank $k$ on $F$ and $Q$ the universal quotient
bundle of rank $n+1-k$.  As is well-known, the relative cotangent
bundle $T^*_{F/S}$ is isomorphic to $\Hom_F(Q, T)$. Thus it is
natural to construct deformations of $T^*_{F/S}$ inside the
endomorphism bundle $\End_F \psi^* V = \psi^* \End_S V$.

Consider the vector bundle $\E$ over $F$ defined as follows: For
$x \in F$,
\begin{equation*}
\E_x:=\{(p, t) \in \End V_{\psi(x)} \times \cit \mid \im\,p
\subset T_x,\ p|_{T_x} = t \Id_{T_x} \}.
\end{equation*}

We have an inclusion $T^*_{F/S, x} = \Hom(Q_x, T_x) \to \E_x$
which sends $q \in \Hom(Q_x, T_x)$ to $(\tilde{q}, 0) \in \E_x$,
where $\tilde{q}$ is the composition
$$
\tilde q: V_{\psi(x)} \to Q_x \xrightarrow{q} T_x \hookrightarrow
V_{\psi(x)}.
$$
The projection to the second factor $\pi: \E \to \cit$ is then an
one-dimensional deformation of $\pi^{-1}(0) = T^*_{F/S}$.

Equivalently, the Euler sequence $0 \to T \to \psi^* V \to Q \to
0$ leads to
$$
0 \to T^*_{F/S} = \Hom_F(Q, T) \to \Hom_F(\psi^*V, T) \to \End_F T
\to 0.
$$
The deformation is simply the inverse image of $\mathbb{C}\Id_F T
\cong \mathbb{C}$.

The projection to the first factor, followed by $\psi_*$
$$
\xymatrix{\E \ar[rd] \ar[r] & \psi^* \End_S V \ar[r]^{\psi_*}
\ar[d] & \End_S V \ar[d] \\ & F \ar[r]^\psi & S}
$$
gives rise to a map $\E \to \End_S V$, which is a birational
morphism onto its image $\bar \E$. Indeed, $\Psi: \E \to \bar \E$
is isomorphic over the loci with ${\rm rank}\,p = k$. In
particular, it is isomorphic outside $\pi^{-1}(0)$. For any $s \in
S$,
$$
\bar \E_s:=\overline{ \{\,p \in \End V_s \mid \mbox{${\rm rank}\,p
= k$ and $p^2 = t p$ for some $t \in \cit$} \,\}}
$$
is the {\em cone of scaled projectors} with rank at most $k$. Thus
$\pi = \bar \pi \circ \Psi$, where
$$
\bar \pi: \bar \E \to \cit \quad \mbox{via}\quad \phi \mapsto
\frac{1}{k}\, {\rm Tr}\, \phi.
$$

For $t \neq 0$, $\Psi_t: \E_t \mathop{\to}\limits^\sim \bar \E_t$.
For $t = 0$, $\psi := \Psi_0: T^*_{F/S} = \E_0 \to \bar\E_0$ is
the open local model of an $A_{n, k}$ contraction.

We do a similar construction for the dual bundle $V^* \to S$.
Under the canonical isomorphism $\End_S V \simeq \End_S V^*$, we
see that $\bar \E$ is identified with $\bar \E' = \bar \E(V^*)$.
Thus we get a birational map $\F: \E \dasharrow \E'$ over $\bar
\E$. This proves

\begin{proposition}
The birational map $\F$ over $\mathbb{C}$:
$$
\xymatrix{\E \ar[rd]^\Psi \ar[rdd]_\pi \ar@{.>}[rr]^\F & & \E'
\ar[ld]_{\Psi'} \ar[ldd]^{\pi'}\\ & \bar \E \ar[d] & \\ &
\mathbb{C} &}
$$
deforms the birational map ($A_{n, k}$ flop) $f: T^*_{F/S}
\dasharrow T^*_{F'/S}$ into isomorphisms.
\end{proposition}

Let $\Gamma$ be the graph closure of $\F: \E \dasharrow \E'$ and
$\Gamma_0$ be its central fiber. By Proposition \ref{equi}, the
map $\Gamma_*: A^*(\E) \to A^*(\E')$ is an isomorphism. Since
$\Gamma \to \E \times_{\bar \E} \E'$ is birational, $(\E
\times_{\bar \E} \E')_*: A^*(\E) \to A^*(\E')$ is again an
isomorphism. It follows that its central fiber $\T_{open}:=
T^*_{F/S} \times_{\bar \E_0} T^*_{F '/S}$ induces an isomorphism
$A^*(T^*_{F/S}) \to A^*(T^*_{F'/S})$.

Consider the fiber product
$$
\T_{loc} := \pit_F(T^*_{F/S} \oplus \O) \times_{\pit_S (\bar \E_0
\times \cit)} \pit_{F'}(T^*_{F '/S} \oplus \O )
$$
and $\T_{\infty}:= \pit_F(T^*_{F/S}) \times_{\pit_S (\bar \E_0)}
\pit_{F'}(T^*_{F'/S}).$ Note that the push-forward map
$A_*(\pit_F(T^*_{F/S})) \to A_*(\pit_F(T^*_{F/S} \oplus \O))$ is
injective.

\begin{proposition} \label{chow}
We have the following commutative diagrams with exact horizontal
rows (induced by the localization formula in Chow groups): \small
$$
\begin{CD}
0 @>>> A_*(\pit(T^*_{F/S})) @>>> A_*(\pit(T^*_{F/S} \oplus \O))
@>>>  A_*(T^*_{F/S}) @>>>0 \\
@.  @V\T_{\infty}VV  @V\T_{loc} VV  @V\T_{open}VV  @.  \\
0 @>>> A_*(\pit(T^*_{F'/S})) @>>> A_*(\pit(T^*_{F'/S} \oplus \O))
@>>>  A_*(T^*_{F'/S}) @>>>0 \\
\end{CD}
$$
\normalsize Thus if $\T_{\infty}$ is an isomorphism, so is
$\T_{loc}$.
\end{proposition}

Note that $\pit_F(T^*_{F/S}) \dasharrow \pit_{F'}(T^*_{F'/S})$ is
a stratified Mukai flop of type $A_{n-2, k-1}$. This allows us to
perform inductive argument later.

\subsection{Deformations of projective local models}
Consider the rational map $\pi: \pit(\E \oplus \O) \dasharrow
\pit^1$ which extends the map $\pi: \E \to \cit$ and maps
$\pit(\E)$ to $\infty$. The map $\pi$ is undefined exactly along
$E:=\pit(\E_0)$. Blow up $\pit(\E \oplus \O)$ along $E$ resolves
the map $\pi$, thus we obtain $\hat \pi: \X := {\rm Bl}_E \pit(\E
\oplus \O) \to \pit^1$.

Since $E$ is a divisor of the central fiber $\pit(\E \oplus \O)_0
= \pit(T^*_{F/S} \oplus \O)$, we have $\hat\pi^{-1}(0) \simeq
\pit(T^*_{F/S} \oplus \O)$. When $t \neq 0$, $\X_t \simeq
\pit(\E)$ which is the compactification of $\E_t$ by $E \cong
\pit(\E_0) = \pit(T^*_{F/S})$. This gives a deformation of
$\pit(T^*_{F/S} \oplus \O)$ over $\pit^1$ with other fibers
isomorphic to $\pit(\E)$.

We do a similar construction on the dual side, which gives a
deformation of $\pit(T^*_{F'/S} \oplus \O)$ by $\hat \pi': \X' \to
\pit^1$. We get also an induced birational map $\F: \X \dasharrow
\X'$ over $\bar\X$ extending $\F: \E \dasharrow \E'$ over $\bar
\E$.

The flop $\F_t: \X_t \dasharrow \X'_t$ for $t \ne 0$ has the
property that there are smooth divisors $E \subset \X_t$ and $E'
\subset \X'_t$ such that (i) the exceptional loci $Z \subset \X_t$
(resp.~$Z' \subset \X'_t$) is contained in $E$ (resp.~$E'$), (ii)
$\F_t|_E: E \dasharrow E'$ is a stratified Mukai flop. We call
such flops \emph{stratified ordinary flops} if furthermore (iii)
$N_{E/\X_t}|_{\mathbb{P}^1} \cong \mathscr{O}$ and
$N_{E'/\X'_t}|_{\mathbb{P}^1} \cong \mathscr{O}$ along the
flopping extremal rays.

If $\F_t|_E$ is of type $A$, $D$ or $E$, then we say $\F_t$ is of
type $A^*$, $D^*$ or $E^*$ respectively. Notice that stratified
ordinary flops of type $A^*_{m, 1}$ are precisely ordinary
$\mathbb{P}^m$ flops, which explains the choice of terminology.

\begin{proposition}\label{proj-deform}
The birational map $\F$ over $\pit^1$:
$$
\xymatrix{\X \ar[rd]^{\Psi} \ar[rdd]_{\hat\pi} \ar@{.>}[rr]^\F
& & \X' \ar[ld]_{\Psi'} \ar[ldd]^{\hat\pi'}\\ & \bar\X \ar[d] & \\
& \mathbb{P}^1 &}
$$
deforms the $A_{n, k}$ flop $f = \F_0: \pit(T^*_{F/S} \oplus \O)
\dasharrow \pit(T^*_{F'/S} \oplus \O)$ into $A^*_{n - 2, k - 1}$
flops $\F_t: \pit(\E) \dasharrow \pit(\E')$ for $t \ne 0$.
\end{proposition}

\begin{proof}
It remains to check condition (iii), which is equivalent to $(E.C)
= 0$ for any flopping curve $C \cong \mathbb{P}^1$. Since $(E.C)$
is independent of $t \in \mathbb{P}^1$ we may compute it at $t =
0$. As a projective bundle $\rho: \X_0 \to F$ it is clear that
$$
K_{\X_0} = -(2\dim F/S + 2)E + \rho^* K_{\X_0}|_F.
$$
Since $(K_{\X_0}.C) = 0$ by the definition of flops, we get $(E.C)
= 0$ as well.
\end{proof}

Clearly for $t \neq 0$, the map $\F_t: \pit(\E) \dasharrow
\pit(\E')$ is an isomorphism only for the case of ordinary Mukai
flop, i.e.\ $k = 1$.
\begin{remark} \label{non-deformable}
 For $A_{n, 2}$ flops, the
deformed flop $\F_t$ is a family of ordinary flop, which has
defect of cup product by \cite{LLW}. As the classical cohomology
ring is invariant under deformations, the fiber product of $f$
does not preserve the ring structure. This implies that we cannot
deform a projective local stratified Mukai flop of type $A_{n, 2}$
into isomorphisms, which is a crucial difference to usual Mukai
flops. Thus there exist quantum corrections even in this local
case.
\end{remark}

\begin{corollary}\label{chow-local}
For projective local model of $A_{n, k}$ flops $X
\xrightarrow{\psi} \bar X \xleftarrow{\psi'} X'$, the
correspondence defined by fiber product $\T = X \times_{\bar X}
X'$ induces isomorphism of Chow motives $[X] \cong [X']$ which
preserves also the Poincar\'e pairing.
\end{corollary}

\begin{proof}
By definition, the $A_{n, k}$ contraction satisfies
$$
2 \codim D = \codim \psi(D)
$$
for each irreducible component $D$ of $E_i$. The deformation
$$
\X \xrightarrow{\Psi} \bar \X \xleftarrow{\Psi'} \X'
$$
constructed by Proposition \ref{proj-deform} is not isomorphic on
general fibers, instead it gives $A^*_{n - 2, k -1}$ flops. Thus
the additional deformation dimension makes it satisfying the
assumption of Proposition \ref{equi}. The result follows by
noticing that the graph closure restricts to $\T$ on the central
fiber.
\end{proof}

\begin{remark}
Proposition \ref{proj-deform} suggests certain inductive structure
on $A_{n, k}$ flops. It will become more useful (e.g.\ for the
discussion of global $A_{n, k}$ flops or Gromov-Witten theory)
after we develop detailed analysis on correspondences.
\end{remark}

\section{Degeneration of correspondences}

\subsection{Setup of degeneration}

Let $f: X \dashrightarrow X'$ be a stratified Mukai flop, say of
type $A_{n, k}$ with $2k < n + 1$. The aim of the following
theorem is to show that the degeneration to normal cone for $(X,
F_k)$ and $(X', F_k')$ splits the correspondence $\T^f$ defined by
$X\times_{\bar X} X'$ into the one $\T^g$ defined by
$Y\times_{\bar Y} Y'$ of type $A_{n - 2, k - 1}$ and its version
$\T^f_{\rm loc}$ on projectivized local models relative to $F_k$
and $F_{k'}$. Conversely we may define $\T$ inductively by gluing
these two parts. Here is the blow-up diagram
$$
\xymatrix{Y = {\rm Bl}_{F_k} X \ar[d]_\phi \ar@{.>}[rr]^g
\ar[rd]& & Y'= {\rm Bl}_{F_k'} X' \ar[d]^{\phi'} \ar[ld] \\
X\ar@{.>}[rr]|{\qquad\qquad\qquad\qquad} \ar[rd]^\psi
& \bar Y = {\rm Bl}_{S_k} \bar X \ar[d] & X' \ar[ld]_{\psi'} \\
& \bar X &}
$$
with $g:Y \dashrightarrow Y'$ being the induced $A_{n - 2, k - 1}$
flop. To save the notations we use the same symbol $\T$ for $\T^f$
and its local models as well if no confusion is likely to arise.

We consider degenerations to the normal cone $W \to \A^1$ of $X$,
where $W$ is the blow-up of $X \times \A^1$ along $F_k \times
\{0\}$. Similarly we get  $W' \to \A^1$ for $X'$.

Note that the central fiber
$$
W_0 := Y_1 \cup Y_2 = Y \cup X_{loc}, \qquad W'_0 := Y_1' \cup
Y_2' = Y' \cup X'_{loc},
$$
where $X_{loc}= \pit (T^*_{F_k/S} \oplus \O)$ and $X'_{loc}= \pit
(T^*_{F'_k/S} \oplus \O)$. The intersections $E:=Y \cap X_{loc} $
and $E':= Y' \cap X'_{loc}$ are isomorphic respectively to
$\pit(T^*_{F_k/S})$ and $\pit(T^*_{F'_k/S})$. The map $f: X
\dasharrow X'$ induces the Mukai flop of the same type for local
models: $f: X_{loc} \dasharrow X'_{loc}$ and Mukai flop $g: Y
\dasharrow Y'$ of type $A_{n - 2, k - 1}$. Let $p: X_{loc} \to
F_k$ be the projection and similarly we get $p'$.

\subsection{Correspondences}

A  {\em lifting} of an element $a \in A^* (X)$ is a couple $(a_1,
a_2)$ with $a_1 \in A^*(Y)$ and $a_2 \in  A^*(X_{loc})$ such that
$\phi_* a_1 + p_* a_2 = a$ and $a_1|_E=a_2|_E$. Similarly one
defines the lifting of an element in $A^*(X')$.

\begin{theorem}\label{coh-red}
Let $a \in A^*(X)$ with $(a_1, a_2)$ and $(a_1', a_2')$ being
liftings of $a$ and $\T a$ respectively. Then
$$
\T a_1 = a_1' \Longleftrightarrow \T a_2 = a_2'.
$$
Moreover it is always possible to pick such liftings.
\end{theorem}

It is instructive to re-examine the Mukai case ($k = 1$) first. In
this case $Y = Y'$ and $f$ is an isomorphism outside the blow-up
loci $Z = F_1$ and $Z' = F_1'$. Let us denote $\T = \T_0 + \T_1$
with $\T_0 = [\bar \Gamma_f] = \phi'_* \phi^*$ and $\T_1$ the
degenerate correspondence $[Z\times_{S} Z']$.

By Lemma 4.2 in \cite{LLW}, it is enough to prove the result for
any single choice of $a_1 = a_1'$. Consider the standard liftings
$$
a(0) = (a_1, a_2) = (\phi^*a, p^*(a|_Z)),
$$
$$
(\T a)(0) = (\phi'^*\T a, p'^*(\T a|_{Z'})).
$$
Since $\phi'^*\T a = \phi^*a + \lambda$ with $\lambda$ supported
on $E' = E$, we may select lifting $(a_1', a_2')$ with $a_1' =
a_1$. In that case,
$$
(\T a_2 - a_2')|_{E'} = a_2|_E - a_2'|_{E'} = a_1|_E - a_1'|_E = 0
$$
by the compatibility constraint on $E$ and the fact that $\T$
restricts to an isomorphism on $E$. Hence $\T a_2 - a_2' = \iota_*
(z')$, for some $z' \in A^*(Z')$, where $\iota: Z' \to X'_{loc}$
is the natural inclusion.

To prove that $z' = 0$, consider
$$
z' = p'_* \iota_* z' = p'_* \T a_2 - p'_* a_2'.
$$
By substituting $\phi'_*a_1' + p'_*a_2' = \T a$, $a_1' = \phi^*a$
and $a_2 = p^*(a|_Z)$, we get
$$
p'_* a_2' = \T a - \phi'_* \phi^* a = \T_1 a.
$$

Let $q, q' $ be the projections of $Z \times_S Z'$ to the two
factors and
$$
j: X_{loc} \times_{\bar{X}_{loc}} X'_{loc} \to Z \times_S Z'
$$
the natural morphism. Then
$$
z' = p'_*\T_{loc} p^*(a|_Z) - \T_1 a = q'_* j_* j^* q^*(a|_Z) - q'
_*q^*(a|_Z).
$$

Note that there exists a unique irreducible  component in $X_{loc}
\times_{\bar{X}_{loc}} X'_{loc}$ birational to $Z \times_S Z'$ via
$j$, so $j_*j^* = {\rm Id}$, which gives $z' = 0$.
\bigskip

Now we proceed for general $A_{n, k}$ flops. It is enough to prove
the result for any single choice of $a_1$ and $a_1'$, since other
choices differ from this one by elements supported on $E$ and $E'$
where the theorem holds by induction on $k$. To make $a_1' = \T
a_1$, notice that $g$ is an isomorphism outside $\tilde F_1 =
\phi_*^{-1}(F_1)$ and $\tilde F_1' = \phi_*'^{-1}(F_1')$ but we
may adjust the standard lifting $\phi'^* \T^f\, a$ only by
elements lying over $F_k'$, namely classes in $E' =
\mathbb{P}(T^*_{F_k'/S_k})$.

The following simple observation resolves this as well as later
difficulties. Recall that $\T = \sum_j \T_j$ with $\T_j = [F_j
\times_{S_j} F_j']$.

\begin{lemma}
We have decomposition of correspondences:
$$
\T^f\, = \phi'_*\T^g\, \phi^* + \T^f_k.
$$
In particular, $\phi'^* \T^f\, = \T^g\, \phi^*$ modulo $A^*(E')$.
\end{lemma}

\begin{proof}
This follows from the definition and the base change property of
fiber product.
\end{proof}

Thus we may pick
$$
a_1 = \phi^* a, \qquad a_1' = \T^g\, (\phi^* a) = \T a_1.
$$
Then
\begin{align*}
(\T a_2 - a_2')|_{E'} &= \T(a_2|_E) - a_2'|_{E'} \\
&= \T(a_1|_E) - a_1'|_{E'} = (\T a_1 - a_1')|_{E'} = 0
\end{align*}
and so $\T a_2 - a_2' = \iota_* (z')$, for some $z' \in
A^*(F_k')$, where $\iota: Z' \to X'_{loc}$ is the natural
inclusion.

To prove that $z' = 0$, consider
$$
z' = p'_* \iota_* z' = p'_* \T a_2 - p'_* a_2'.
$$
By substituting $\phi'_*a_1' + p'_*a_2' = \T a$, $a_1' =
\T\phi^*a$ and $a_2 = p^*(a|_{F_k})$, we get
$$
p'_* a_2' = \T a - \phi'_* \T \phi^* a = \T_k a
$$
by the above lemma.

Let $q, q' $ be the projections of $F_k \times_{S_k} F_k'$ to the
two factors and
$$
j: X_{loc} \times_{\bar{X}_{loc}} X'_{loc} \to F_k \times_{S_k}
F_k'
$$
the natural morphism. Then
$$
z' = p'_*\T_{loc} p^*(a|_{F_k}) - \T_k a = q'_* j_* j^*
q^*(a|_{F_k}) - q' _*q^*(a|_{F_k}).
$$

Note that there exists a unique irreducible  component in $X_{loc}
\times_{\bar{X}_{loc}} X'_{loc}$ birational to $F_k \times_{S_k}
F_k'$ via $j$, so $j_*j^* = {\rm Id}$, which gives $z' = 0$. The
proof is complete.

\subsection{Cup product and the Poincar\'e pairing}

Besides correspondences, we also need to understand the effect on
the Poincar\'e pairing under degeneration. We will in fact
degenerate classical cup product and this works for any
degenerations to normal cones $W \to X \times \mathbb{A}^1$ with
respect to $Z \subset X$. Let $W_0 = Y_1 \cup Y_2$, where $\phi:
Y_1 = Y \to X$ is the blow up along $Z$, $p: Y_2 = \tilde E =
\mathbb{P}_Z(N_{Z/X} \oplus \O) \to Z$ is the local model and $Y_1
\cap Y_2 = E$ is the $\phi$ exceptional divisor. Let $i_1: E
\hookrightarrow Y_1$, $i_2: E \hookrightarrow Y_2$.

\begin{lemma}\label{poincare}
Let $a, b \in H^*(X)$. Then for any lifting $(a_1, a_2)$ of $a$
and any lifting $(b_1, b_2)$ of $b$, the pair $(a_1b_1, a_2b_2)$
is a lifting of $ab$.

In particular, if $a$, $b$ are of complementary degree, then we
have an orthogonal splitting of the Poincar\'e pairing: $(a.b)_X =
(a_1.b_1)_{Y_1} + (a_2.b_2)_{Y_2}$.
\end{lemma}

\begin{proof}
We compute
$$
a.b = \phi_* a_1.b + p_*a_2.b = \phi_*(a_1.\phi^*b) +
p_*(a_2.p^*(b|_Z)).
$$

Since $a_1.\phi^*b|_E = a_2.p^*(b|_Z)|_E$, $(a_1 b_1, a_2 b_2)$ is
a lifting of $ab$ for the special lifting $(b_1, b_2) = (\phi^* b,
p^*(b|_Z))$ of $b$. By \cite{LLW}, Lemma 4.2, any other lifting of
$b$ is of the form $(b_1 + i_{1*}e, b_2 - i_{2*}e)$ for some class
$e$ in $E$.

Since $i_1^*a_1.e = i_2^*a_2.e$ is a class $e' \in H^*(E)$. The
correction terms are
$$
a_1.i_{1*}e = i_{1*}(i_1^* a_1.e) = i_{1*}e', \quad -a_2.i_{2*}e =
-i_{2*}(i_2^* a_2.e) = -i_{2*}e'.
$$

The lemma then follows from
$$
i_1^* (i_{1*} e') = e'.c_1(N_{E/Y_1}) = -e'.c_1(N_{E/Y_2}) = i_2^*
(-i_{2*} e')
$$
and $\phi_* i_{1*}e' - p_* i_{2*}e' = (\phi|_E)_* e' - (\phi|_E)_*
e' = 0$.
\end{proof}

\begin{theorem}[{\bf $=$ Theorem \ref{main}}]\label{chow-main}
For $A_{n, k}$ flops, $\T$ induces an isomorphism on Chow motives
and the Poincar\'e pairing.
\end{theorem}

\begin{proof}
If $f: X \dashrightarrow X'$ is an $A_{n, k}$ flop, $f\times {\rm
id}: X \times T \dashrightarrow X' \times T$ is also an $A_{n, k}$
flop. Thus by the identity principle, to prove $[X] \cong [X']$ we
only need to prove the equivalence of Chow groups under $\T$ for
any $A_{n, k}$ flop.

We prove this for all $n$ with $n + 1 > 2k$ by induction on $k$.
We start with $k = 0$, which is trivial.

Given $k \ge 1$, by Theorem \ref{coh-red}, the equivalence of Chow
groups is reduced to the $A_{n - 2, k - 1}$ case and the local
$A_{n, k}$ case. The former is true by induction. The later
follows from Corollary \ref{chow-local} directly.

The same procedure proves the isomorphism of Poincare pairings by
using Theorem \ref{coh-red}, Lemma \ref{poincare} and Corollary
\ref{chow-local}.
\end{proof}

For cohomology rings we need to proceed carefully. In order to run
induction on $k$, using Theorem \ref{coh-red}, Lemma
\ref{poincare} we must first consider the local $A_{n, k}$ case.
By remark \ref{non-deformable}, for $k=2$, the classical cup
product is \emph{not} preserved by the correspondence $\T$! This
is analyzed in the next section.

\section{Quantum corrections}

\subsection{The proof of Theorem \ref{main-2}}

We now prove the invariance of big quantum product attached to the
extremal rays, up to analytic continuations, under $A_{n, 2}$
flops.

As in the precedent section, we consider degenerations to the
normal cone $W \to \A^1$ of $X$ and  $W' \to \A^1$ of $X'$. Note
that the map $f: X \dasharrow X'$ induces the Mukai flop of the
same type for local models: $f: X_{loc} \dasharrow X'_{loc}$ and
Mukai flop $g: Y \dasharrow Y'$ of type $A_{n - 2, 1}$.

By the degeneration formula (for the algebraic version used here,
c.f.\ \cite{Li}), any Gromov-Witten invariant on $X$ splits into
products of relative invariants of $(Y, E)$ and $(X_{loc}, E)$.
Let $a \in H^*(X)^{\oplus m}$ with lifting $(a_1, a_2)$:
$$
\langle a \rangle^X_{g, m, \beta} = \sum_{\eta \in \Omega(g,
\beta)} C_\eta \left(\langle a_1 \rangle^{\bullet (Y,
E)}_{\Gamma_1}. \langle a_2 \rangle^{\bullet (X_{loc},
E)}_{\Gamma_2}\right)^{E^\rho}.
$$
Here $\rho$ is the number of gluing points (in $E$) and $\Gamma_1
\cup \Gamma_2$ forms a connected graph. Thus $\rho = 0$ if and
only if one of the $\Gamma_i$ is empty.

The relative invariants take values in $H^*(E^\rho)$ and the
formula is in terms of the Poincar\'e pairing of $E^\rho$.

We apply it to $X'$ as well and get:
$$
\langle \T a \rangle^{X'}_{g, m, \T\beta} = \sum_{\eta' \in
\Omega(g, \T\beta)} C_{\eta'} \left(\langle a_1' \rangle^{\bullet
(Y', E')}_{\Gamma_1'}. \langle a_2' \rangle^{\bullet (X'_{loc},
E')}_{\Gamma_2'}\right)^{E'^\rho}.
$$

There is a one to one correspondence between admissible triples
$\eta = (\Gamma_1, \Gamma_2, I_\rho)$ and $\eta' = (\Gamma_1',
\Gamma_2', I_{\rho'})$ via $\eta' := \T \eta$. The combinatorial
structure is kept the same, while the curve classes are related by
$\T$. We do still need the cohomology class splitting on $X$ and
$X'$ be to compatible.

By Theorem \ref{coh-red} we may split the cohomology classes $a
\in H^*(X)^{\oplus m}$ into $(a_1, a_2)$ with $a_i \in
H^*(Y_i)^{\oplus m}$ and $\T a \in H^*(X')^{\oplus m}$ into
$(a_1', a_2')$ with $a_i' \in H^*(Y_i')^{\oplus m}$, such that
$$
\T a_1 = a_1',\qquad \T a_2 = a_2'.
$$

By Theorem \ref{chow-main}, the Poincar\'e pairing is preserved by
$\T$ under stratified Mukai flops $E \dasharrow E'$, the same
holds true for $E^\rho \dasharrow E'^\rho$ by $\T^\rho$, which for
simplicity is still denoted by $\T$. Thus by the degeneration
formula the problem is reduced to showing that $\T$ maps the
relative invariants of $(Y, E)$ and $(X_{loc}, E)$ to the
corresponding ones of $(Y', E')$ and $(X'_{loc}, E')$.

Since we are only interested in invariants attached to the
extremal ray $\beta = d\ell$, for any splitting $\beta = (\beta_1,
\beta_2) = (d_1\ell, d_2\ell)$, we must have $\rho = (E.\beta_2) =
d_2 (E.\ell) = 0$ (since $\ell$ can be represented by a curve in
$F_2$). But this implies that $\beta$ is not split at all and in
the degeneration formula the invariant $\langle a \rangle^{X}_{g,
m, d\ell}$ goes to $Y$ or $X_{loc}$ completely:
$$
\langle a \rangle^{X}_{g, m, d\ell} = \langle a_1 \rangle^{Y}_{g,
m, d\ell} + \langle a_2 \rangle^{X_{loc}}_{g, m, d\ell}.
$$

\begin{lemma}\label{abs-rel}
$\T$ maps isomorphically the cup product and full Gromov-Witten
theory of $Y$ to those of $Y'$. Moreover $\langle a_1
\rangle^{Y}_{g, m, d\ell} = 0$ for all $d \in \mathbb{N}$.
\end{lemma}

\begin{proof}
The birational map $g: Y \dasharrow Y'$ is a Mukai flop of type
$A_{n-2, 1}$. Hence this follows from \cite{LLW}, Theorem 6.3.

Indeed this follows from Lemma \ref{poincare} and the above
degeneration formula by applying it to the Mukai flop $Y
\dasharrow Y'$. Here we use the facts that projective local models
of Mukai flops $g_{loc}: Y_{loc} \dasharrow Y'_{loc}$ can be
deformed into isomorphisms $g_t: Y_t \mathop{\to}\limits^\sim
Y'_t$ and that the cup product as well as the Gromov-Witten theory
are both invariant under deformations. For $\ell$ being the
extremal ray of $g_{loc}$, if $d\ell \sim C_t$ for $t \ne 0$ then
$C'_t \cong g_t(C_t) \sim \T d\ell = -d\ell'$, which is
impossible. Thus $\langle a_1 \rangle^{Y}_{g, m, d\ell} = 0$ for
all $d \in \mathbb{N}$.
\end{proof}

Denote by $\langle a \rangle_f = \sum_{d = 0}^\infty \langle a
\rangle_{0, m, d\ell}\, q^{d\ell}$, the generating function of $g
= 0$ Gromov-Witten invariants attached to the extremal ray. Then
the degeneration formula and Lemma \ref{abs-rel} lead to
$$
\langle a \rangle^X_f = \delta_{n3}\langle a_1 \rangle^Y_{0,3,0} +
\langle a_2\rangle^{X_{loc}}_f.
$$

The correspondence $\T$ acts on $q^\beta$ by $\T q^\beta = q^{\T
\beta}$. In particular for the extremal rays $\ell$ and $\ell'$ we
have $\T q^{d\ell} = q^{-d\ell'}$. If we regard $q^{\ell} =
e^{-2\pi(\omega.\ell)}$ as an analytic function on $\omega \in
H^{1,1}_\mathbb{R}(X)$, then it is known that $\langle a
\rangle^X_f$ converges in the K\"ahler cone $\mathcal{K}_X$ of
$X$. Under the identification $H^{1,1}_\mathbb{R}(X) \cong \T
H^{1,1}_\mathbb{R}(X) = H^{1,1}_\mathbb{R}(X')$, it makes sense to
compare $\T\langle a \rangle^X_f$ with $\langle \T a
\rangle^{X'}_f$ as analytic functions on $\mathcal{K}_X \cup
\mathcal{K}_{X'} \subset H^{1,1}_\mathbb{R}$ up to analytic
continuations.

\begin{lemma} \label{analytic}
For $m \ge 3$, $\T\langle a_2 \rangle^{X_{loc}}_f \cong \langle \T
a_2 \rangle^{X'_{loc}}_f$ up to analytic continuations.
\end{lemma}

\begin{proof}
By Proposition \ref{proj-deform}, the $A_{n, 2}$ flop $f : X_{loc}
\dasharrow X'_{loc}$ can be projectively deformed into ordinary
$\mathbb{P}^{n - 2}$ flops. By the deformation invariance of
Gromov-Witten theory, the lemma is reduced to the case of ordinary
flops (with non-trivial base). For simple ordinary flops the
invariance
$$
\T\langle a_2 \rangle^{X_{loc}}_f \cong \langle \T a_2
\rangle^{X'_{loc}}_f
$$
up to analytic continuations is proved in \cite{LLW}. It has been
extended to general ordinary flops with base in \cite{LLW2}. Hence
the lemma follows.
\end{proof}

Notice that the $g = 0$, $d = 0$ invariants are non-zero if and
only if $m = 3$ and they are given by the cubic product. By Lemma
\ref{abs-rel}, $\langle a_{11}, a_{12}, a_{13} \rangle^Y_{0,3,0} =
\langle \T a_{11}, \T a_{12}, \T a_{13} \rangle^{Y'}_{0,3,0}$.
From Lemma \ref{analytic} we get $\T \langle a \rangle^X_f =
\langle \T a \rangle^{X'}_f$ for $m \ge 3$.

Together with the $\T$ invariance of Poncar\'e pairing, the big
quantum product attached to the extremal ray is invariant under
$\T$. This completes the proof of Theorem \ref{main-2} for type
$A_{n,2}$. The cases of type $D_5$ and $E_{6, I}$ are completely
similar, since the geometric picture in Proposition
\ref{proj-deform} is the same by \cite{CF}. The proof is complete.

\begin{remark}
The degeneration formula is in terms of the Poincar\'e pairing of
relative GW invariants. Thus \emph{invariance of the Poincar\'e
pairing} is crucial in our study. Indeed, the Poincar\'e pairing
together with 3-point functions determine the (small) quantum
product. So far this is \emph{the only constraint} we have found
for the correspondence $\T$ under $K$ equivalence to be
\emph{canonical}.
\end{remark}

\subsection{A new topological constraint}

Consider a stratified Mukai flops of type $A_{n, 2}$,  $f: X
\dasharrow X'$ with $i: F_2 \hookrightarrow X$, such that $\T$
preserves the cup product (e.g.~for $X$ and $X'$ being
hyper-K\"ahler manifolds). By Proposition \ref{proj-deform} and
Remark \ref{non-deformable}, there exists defect of cup product on
$X_{loc}$. A priori there seems to be a contradiction. A closer
look at them leads to

\begin{proposition} \label{main-3}
For a Mukai flop $f: X \dasharrow X'$ of type $A_{n, 2}$, $D_{5}$
or $E_{6, I}$, if the restriction map $i^*: H^*(X, \mathbb{Q}) \to
H^*(F_2, \mathbb{Q})$ is surjective then $\T$ does not preserve
the cup product. In particular, if $X$ is hyper-K\"ahler then
$i^*$ is not surjective.
\end{proposition}

\begin{proof}
We shall investigate the degeneration analysis on cup product for
an arbitrary $A_{n, 2}$ flop $f$ as presented above. The other
cases are similar.

Let $a=(a_1, a_2)$ and $b=(b_1, b_2)$ be two elements in $H^*(X)$
with their lifting. By Lemma \ref{poincare}, $ab = (a_1b_1,
a_2b_2)$. Then Theorem \ref{coh-red} implies that $\T(ab) =
(\T(a_1b_1), \T(a_2b_2))$. By Lemma \ref{poincare} again
$$
\T(a)\T(b) = (\T(a_1)\T(b_1), \T(a_2)\T(b_2)) = (\T(a_1b_1),
\T(a_2)\T(b_2)),
$$
where the last equality follows from Lemma \ref{abs-rel} applied
to $g: Y \dasharrow Y'$. So
$$
\T(ab) = \T(a)\T(b) \quad \Longleftrightarrow \quad \T(a_2b_2) =
\T(a_2) \T(b_2).
$$

That is, the invariance of cup product on $H^*(X)$ is equivalent
to the invariance on elements in $H^*(X_{loc})$ which come from
lifting of elements in $H^*(X)$. Indeed let $i :F_2
\hookrightarrow X$ and $p: X_{loc} \to F_2$ being the projection,
we may choose standard lifting $a_2 = p^* i^* a$. Such elements
form a subring
$$
\Delta_f := p^*i^* H^*(X) \subset p^* H^*(F_2) \subset
H^*(X_{loc}).
$$

By applying this analysis to the case $X = X_{loc} =
\mathbb{P}_{F_2}(T^*_{F_2/S} \oplus \mathscr{O})$ where the cup
product is not preserved under $\T$, we find that the defect of
cup product is completely realized in the subring $p^* H^*(F_2)$
since
$$
i^* H^*(X_{loc}) = H^*(F_2).
$$

For general $X$, if $i^*: H^*(X, \mathbb{Q}) \to H^*(F_2,
\mathbb{Q})$ is surjective, then $\Delta_f \otimes \mathbb{Q} =
p^* H^*(F_2) \otimes \mathbb{Q}$ must contain the defect on
$X_{loc}$, hence the ring structure on $H^*(X)$ is not preserved
under $\T$. This completes the proof.
\end{proof}

\begin{example}
Consider a simple $A_{4,2}$ flop on hyper-K\"ahler manifold $X$ of
dimension 12, where $F_2 = G(2, 5)$ is of dimension 6. The divisor
$c_1(G) = i^*H$ for some $H \in H^2(X, \mathbb{Q})$. But $c_2(G)
\not\in i^* H^4(X, \mathbb{Q})$.
\end{example}

\end{document}